\documentclass[reqno]{amsart}
\usepackage{setspace}
\usepackage{amsmath}
\usepackage{amssymb}
\makeatletter
\def\bc{\begin{center}}
\def\ec{\end{center}}

\def\s2c{\vskip 2cm}
\def\bt{\begin{Theorem}}
\def\et{\end{Theorem}}
\def\bd{\begin{Definition}}
\def\ed{\end{Definition}}
\def\bl{\begin{Lemma}}
\def\el{\end{Lemma}}
\def\bcor{\begin{Corollary}}
\def\ecor{\end{Corollary}}
\def\bpr{\begin{Proposition}}
\def\epr{\end{Proposition}}

\def\mysection{\setcounter{equation}{0}\section}
\newtheorem{Lemma}{Lemma}[section]
\newtheorem{Theorem}[Lemma]{Theorem}
\newtheorem{Proposition}[Lemma]{Propostion}
\newtheorem{Definition}[Lemma]{Definition}
\newtheorem{Corollary}[Lemma]{Corollary}
\newtheorem{theorem}{Theorem}[section]
\newtheorem{lemma}[theorem]{Lemma}
\newtheorem{corollary}[theorem]{Corollary}

\def\elsartstyle{%
    \def\normalsize{\@setfontsize\normalsize\@xiipt{14.5}}
    \def\small{\@setfontsize\small\@xipt{13.6}}
    \let\footnotesize=\small
    \def\large{\@setfontsize\large\@xivpt{18}}
    \def\Large{\@setfontsize\Large\@xviipt{22}}
    \skip\@mpfootins = 18\p@ \@plus 2\p@
    \normalsize
} \@ifundefined{square}{}{} \makeatother

\author{ Duranta Chutia and Rajib Haloi$^*$}
\date{}
\address{Duranta Chutia \newline
 Department of Mathematical Sciences,
 Tezpur University, Sonitpur,
 ASSAM, Pin- 784028, India}
\email{durantachutia123@gmail.com, Phone +913712-275511,  Fax +913712-267006}

\address{Rajib Haloi \newline
 Department of Mathematical Sciences,
 Tezpur University, Sonitpur,
 ASSAM, Pin- 784028, India}
\email{rajib.haloi@gmail.com, Phone +913712-275511,  Fax +913712-267006}

\thanks{$^*$Corresponding author's e-mail:rajib.haloi@gmail.com.}

\begin{document}

\title{Weighted Inequalities for One-sided Fractional  Minimal Function}
\maketitle \vskip .5cm \noindent {\bf Abstract:}
 {In this article, we present weighted norm inequality for a fractional one-sided minimal function. We prove  weighted weak and strong type norm inequalities for the one-sided  minimal function on $\mathbb{R}.$ We construct two weight classes  to ensure the sufficient and necessary part of the weak and strong type estimates. Later, we establish an equivalence relation between these two weight classes. } 

 \vskip .5cm \noindent {\bf AMS Classification (2010):}{ 42B25.\\
Corresponding author: Rajib Haloi,  e-mail: rajib.haloi@gmail.com, Telephone: +91-03712-275511, Fax: +91-03712-267006.}

\mysection{Introduction}
The Hardy-Littlewood maximal function for a locally integrable function $f$ on $\mathbb{R}^n$ is defined as
$$M(f)(x) = \sup_{x \in Q}{\frac{1}{|Q|}}\int_{Q}|f|dy,$$
where the supremum runs over  all those cubes $Q$ having $x$ and whose sides are parallel to the coordinate axes. Maximal function has a significant role in obtaining convergence of some integral operators and also plays a key role in obtaining the Lebesgue differentiation theorem. After the introduction of $A_p$ weights by  Muckenhoupt \cite{MUCK}, the study of weighted estimate for the maximal function has got much more attention. Some new results with some new techniques have been established for the maximal function and we refer some of the literatures \cite{AND, AND1, CHE, FEF,GAR, GRA, GRA1, GRA2, LER, FJM, FJM2, FJM1, MUCK, OM, LAI, SAW, SAW1, SS} for more detail on maximal function and its applications.

It is quite natural to verify how the weighted estimate will behave if we replace supremum by infimum in the construction of maximal function. In this regard, Cruz-Uribe \textit{et al.} \cite{CRUZ} introduce a new concept called minimal function and later it motivates them to characterize the class of functions  for which the reverse H\"older inequality makes sense. The authors\cite{CRUZ}  define the minimal function as
$$m(f)(x) = \inf_{x \in Q}{\frac{1}{|Q|}}\int_{Q}|f|dy,$$
where the infimum is taken over cubes as similar to the case for maximal function.

In defining the minimal function the locally integrable property of $f$ can be removed by simply assuming  $f$ to be a measurable function on $\mathbb{R}^n$ and it turns out to be a main difference between the maximal and minimal function.  Cruz-Uribe \textit{et al.} \cite{CRUZ} also prove the one weight norm inequality for the minimal function on $\mathbb{R}^n.$ Later Cruz-Uribe \textit{et al.} \cite{CRUZ1} establish the two weight problem for minimal function and prove the result in $\mathbb{R}.$ Interestingly they come out with a beautiful result which establishes the equivalence between the weight class for the weak type and the weight class for the strong type norm inequality which is not the case for maximal function.

The one-sided version of the minimal function has been studied by Cruz-Uribe \textit{et al.} \cite{CRUZ2} and they define the one-sided minimal function, $m^+$ for a measurable function $f$ as
$$m^+(f)(x)= \inf_{h >0} \frac{1}{h} \int_{x}^{x+h}|f|.$$
A backward version of the minimal function, $m^-$ is defined as 
$$m^-(f)(x)= \inf_{h >0} \frac{1}{h} \int_{x-h}^{x}|f|.$$
The authors \cite{CRUZ2} prove two weight weak and strong type estimates for the minimal function on $\mathbb{R}$ with the help of two different weight classes and use these for the study of differentiability of the integral. Further they establish an equivalence relationship between these two weight classes. They also characterize the class of functions which satisfy one-sided reverse H\"older inequality in weak sense.

We plan to study weighted inequalities for the fractional one-sided minimal function on $\mathbb{R}.$ This will extend  the existing results for one-sided minimal function. For this, we define the fractional version of one-sided minimal function as follows.

Let $f$ be a measurable function on $\mathbb{R}.$ For $ 0 \leq \mu < \infty,$ we define a fractional one-sided minimal function of order $\mu$ as
$$m_{\mu}^{+}(f)(x)= \inf_{ h > 0} \frac{1}{ h^{1 + \mu}}\int_{x}^{x+h}\lvert f(t) \rvert dt. $$
If we consider the integral from $x-h$ to $x,$ then we get an another version of one-sided minimal function and precisely it is given as 
$$m_{\mu}^{+}(f)(x)= \inf_{ h > 0} \frac{1}{ h^{1 + \mu}}\int_{x-h}^{x}\lvert f(t) \rvert dt. $$

If we put $\mu =0,$ then we get the one-sided minimal function $m^+.$ Thus $m^+$ is a particular case of $m^+_{\mu}.$ Some literature uses forward minimal function for $m_{\mu}^+$ and backward minimal function for $m_{\mu}^{-}.$ As similar to the case of maximal function the results of forward and backward minimal function can be interchangeable with a little modification in the weight class.

In this note, we mainly deal with functions on $\mathbb{R}.$ By a weight, we mean a non-negative function which takes  positive values \textit{a.e.} on its domain. We write $w(E)=\int_{E}w dy,$ where $E$ is a measurable and $w$ is a locally integrable function. We write $C$ to denote a positive constant not necessarily same in all cases. We use the notation $s^{\prime}$ to denote the conjugate exponent of $s( >1).$

We prepare our article as follows. We state some  already known results for the one-sided minimal function in the section 2. In the section 3, we state our main results and their proofs have been given in the section 4.

\mysection{ Preliminaries}\label{Preliminaries}
In this section, we briefly state some results which are already proved in \cite{CRUZ1, CRUZ2} and we recall some of the basic definitions, lemmas and theorems that are used in the remaining part of the article.

To establish the weak $(p,p),~ p> 0$ estimate for the minimal function, the authors \cite{CRUZ1} define a certain weight class, $W_p$ which is similar to the Muckenhoupt weight class, $A_p$ for the case of Hardy-Littlewood maximal function but with $p$ replaced by $-p.$ Precisely the weight class $W_p$ is as follows.
\begin{Definition}
Let $(U, V)$ be a pair of weight on $\mathbb{R}.$ Given $p >0,$ we say $(U, V) \in W_p$ if there exists $C >0$ such that
$$\frac{1}{|I|}\int_{I} U \leq C \Bigg(\frac{1}{|I|} \int_{I} V^{\frac{1}{p+1}}  \Bigg)^{p+1}$$ 
holds for every interval $I$ in $\mathbb{R}.$
\end{Definition}
With the help of the $W_p$ condition, the following result regarding the weak type estimate has been established.
\begin{theorem}\cite{CRUZ1}
Let $(U, V)$ be a pair of weight on $\mathbb{R}.$ For $p, \lambda >0$ the following statements are equivalent.
\begin{itemize}
\item[(i)] The pair $(U,V) \in W_p.$
\item[(ii)] The weak $(p, p)$ holds for the minimal function. \textit{i.e.}
$$ U\Bigg(\bigg \{ x \in \mathbb{R} : m(f)(x) < \frac{1}{\lambda} \bigg\} \Bigg) \leq \frac{C}{{\lambda}^p} \int_{\mathbb{R}} \frac{V}{|f|^p}$$
holds for some constant $ C > 0.$

\end{itemize}
\end{theorem}
As similar to the case of maximal function, we need an another weight class as similar to the Sawyer's testing type condition  to show the two weight strong type inequality for the minimal function and the weight class to show the strong type inequality is defined in \cite{CRUZ1} as follows.
\begin{Definition}
Let $(U, V)$ be a pair of weight on $\mathbb{R}.$ Given $p >0,$ we say $(U, V) \in W_p^*$ if there exists
 $C >0$ such that
$$\int_{I} \frac{U}{{m(\omega/\chi_{I})}^p} \leq C \int_{I} \omega$$ 
holds for every interval $I$ in $\mathbb{R},$ where $\omega = V^{\frac{1}{p+1}}$
\end{Definition}
In the remaining part of the article we use the notation $\omega = V^{\frac{1}{p+1}}$ for $p> 0.$ The strong type estimate $(p, p),~ p>0$ for the minimal function follows from the weight class $W_p^*$ and this result has been proved in \cite{CRUZ1}.
\begin{theorem}\cite{CRUZ1}
Let $(U, V)$ be a pair of weight on $\mathbb{R}.$ For $p >0$ the following statements are equivalent.
\begin{itemize}
\item[(i)] The pair $(U,V) \in W_p^*.$
\item[(ii)] The strong $(p, p)$ holds for the minimal function. \textit{i.e.}
$$ \int_{\mathbb{R}} \frac{U}{(m(f))^p} \leq C\int_{\mathbb{R}} \frac{V}{|f|^p}$$
holds for some constant $ C > 0.$

\end{itemize}
\end{theorem}

The one-sided minimal function has been considered in \cite{CRUZ2}, where the authors define the one-sided version of the weight class $W_p$ and $W_p^*$ denoted by $W_p^+$ and $(W_p^+)^*$ respectively and with the help of  these weight classes they prove the weak and strong type inequalities for the one-sided minimal function.

The weight class, $W_p^+$ corresponding to the weak weighted inequality for $m^+$ is defined in the following way.
\begin{Definition}
Given $p >0,$ we say a pair of weight $(U, V) \in W_p^+$ if there exists a constant $C >0$ such that
$$\frac{1}{|I^-|}\int_{I^-} U \leq C \Bigg(\frac{1}{|I|} \int_{I} V^{\frac{1}{p+1}}  \Bigg)^{p+1}$$ 
holds for each interval $I=[\alpha, \beta]$ in $\mathbb{R}$ with $2 |I^-|= |I|$ and $I^-=[\alpha, \gamma]$
\end{Definition}
Cruz-Uribe \textit{et al.} \cite{CRUZ2} establish the two weight weak $(p, p)$ for $m^+$ with the help of $W_p^+$ weight  and their main result is the following.
\begin{theorem}\cite{CRUZ2}
Let $p, \lambda >0.$ Then the pair $(U,V)\in W_p^+$ if and only if the weak $(p,p)$ inequality
$$ U\Bigg(\bigg \{ x \in \mathbb{R} : m^+(f)(x) < \frac{1}{\lambda} \bigg\} \Bigg) \leq \frac{C}{{\lambda}^p} \int_{\mathbb{R}} \frac{V}{|f|^p}$$
holds for some constant $C >0.$
\end{theorem}
An inequality similar to the one-sided Sawyer's testing type condition is used to show the two weighted strong $(p,p)$ for $m^+.$ The condition is denoted by $(W_p^+)^*$ and it is defined as follows.
\begin{Definition}
Given $p > 0,$ we say the pair $(U, V) \in (W_p^+)^*$ if
$$\int_{I} \frac{U}{m^+(\omega/\chi_I)^p} \leq C \int_{I}\omega$$
holds for some constant $C >0$ and for any interval $I$ in $\mathbb{R}.$
\end{Definition}
The two weighted strong $(p, p)$ for one-sided minimal function is obtained using the $(W_p^+)^*$ condition and the main result is the following.

\begin{theorem}\cite{CRUZ2}
Given $p >0$ and for a pair of weight $(U, V)$ on $\mathbb{R}$ the following statements are equivalent.
\begin{itemize}
\item[(i)] The pair $(U,V) \in (W_p^+)^*.$
\item[(ii)] The strong $(p, p)$ holds for the one-sided minimal function. \textit{i.e.}
$$ \int_{\mathbb{R}} \frac{U}{(m^+(f))^p} \leq C\int_{\mathbb{R}} \frac{V}{|f|^p}$$
holds for some constant $ C > 0.$

\end{itemize}
\end{theorem}


\mysection{ Main Results}\label{Main Results}
In this section, we define two weight classes to prove the weak and strong type weighted estimate for the fractional minimal function and we prove our main results. We use techniques and ideas from \cite{CRUZ2, SAW} in proving our results. We start this section with a basic lemma whose proof can be obtained with a suitable modification of the Lemma $2.1$ in \cite{CRUZ2}.
\begin{lemma}\label{fml1}
Let $I$ be a interval in $\mathbb{R}.$ We consider a collection of intervals $\{I_a\}_{a}$ contained in $I$ such that given a function $w$ and for $0 \leq \mu < \infty,$ $\int_{I_a}wdx \leq C |I_a|^{1 + \mu},$ for each $a.$ If $J =\cup_{a}I_a,$ then $\int_{J}w dx \leq C (2|J|)^{1+\mu}.$
\end{lemma}
Next we define the fractional one-sided weight class, $W_{p, q}^+,$ $0 < p \leq q < \infty,$  for proving the two weight weak $(p, q)$ estimate for $m^+_{\mu}.$ 
\begin{Definition}\label{fmd1}
Given $0 \leq \mu < \infty$ and $0 < p \leq q < \infty,$ we say a pair of weight $(U, V) \in W_{p,q}^{+}$ if 
$$\frac{1}{|I^-|}\int_{I^-} U \leq \frac{C}{|I|^{1 + (\mu - \frac{1}{p})q}} \Bigg(\frac{1}{|I|} \int_{I} V^{\frac{1}{p+1}}  \Bigg)^{\frac{(p+1)q}{p}}$$ 
holds for each interval $I=[\alpha, \beta]$ in $\mathbb{R}$ with $2 |I^-|= |I|$ and $I^-=[\alpha, \gamma].$\end{Definition}
This definition is about two intervals $I^-$  and $I$ with the relation $2|I^-|=|I|$ and the property that their left points are the same. Next we give a more generalized version of the Definition \ref{fmd1}. For this we consider any subinterval $I^-$ of $I$ with the same starting point and with the property $0 < \frac{|I^-|}{|I|} < 1.$ 
\begin{Definition}\label{fmd2}
Let $I^-=[\alpha, \gamma]$ be any subinterval of $I=[\alpha, \beta].$ We set $\eta =\frac{|I^-|}{|I|}.$ Given $0 \leq \mu < \infty$ and $0 < p \leq q < \infty,$ we say a pair of weight $(U, V) \in W_{p,q, \eta}^{+}$ if 
$$\frac{1}{|I^-|}\int_{I^-} U \leq \frac{C}{ \eta(1 - \eta)^{(1 + \mu)q}}\frac{1}{ |I|^{1 + (\mu - \frac{1}{p})q}} \Bigg(\frac{1}{|I|} \int_{I} V^{\frac{1}{p+1}}  \Bigg)^{\frac{(p+1)q}{p}}$$ 
holds for some constant $C >0.$
\end{Definition}
Next we relate these two weight class $W_{p, q}^+$ and $W_{p, q, \eta}^+.$ To do this we need a decomposition of an interval in $\mathbb{R},$ popularly known as plus-minus decomposition. This idea was introduced in \cite{CRUZ2} and it is given in the following way.
\begin{Definition}\label{fmd3}
Let $I = [\alpha, \beta]$ be a finite interval. We define a sequence $\{x_k\}_{k \geq 0}$ recursively from the interval $I$ as, set $x_0=\alpha$ and for each $k \geq 1$  we define
$$x_k=\frac{\beta + x_{k-1}}{2}.$$
For $k \geq 1,$ we construct three subintervals of $I$ from the sequence $\{x_k\}$ as $J_k^-=[x_{k-1}, x_k], J_k^+=[x_k, x_{k+1}]$ and $J_k=[x_{k-1}, x_{k+1}].$

From the construction itself $I =\cup_{k \geq 0}J_k^-.$
\end{Definition}
The next result is about the relationship between $W_{p, q}^+$ and $W_{p, q, \eta}^+$ and the summary of the following theorem says that these two conditions are equivalent.
\begin{theorem}\label{fmt1}
A pair of weight $(U, V) \in W_{p, q}^+$ if and only if $ (U, V) \in W_{p,q,\eta}^+.$
\end{theorem}
Now we are ready to state the two weight weak estimate for the function $m_{\mu}^+.$ Theorem \ref{fmt1} makes our task easier to prove the following result. 
\begin{theorem}\label{fmt2}
Let $ 0 < p \leq q < \infty$ and $0 \leq \mu < \infty.$ Then for $\lambda >0$ the following conditions are equivalent.
\begin{itemize}
\item[(i)]  The pair of weight $(U, V) \in W_{p,q}^+.$
\item[(ii)] The two weighted weak $(p,q)$ inequality
$$ U\Bigg(\bigg \{ x \in \mathbb{R} : m_{\mu}^+(f)(x) < \frac{1}{\lambda} \bigg\} \Bigg) \leq \frac{C}{{\lambda}^q} \Bigg(\int_{\mathbb{R}} \frac{V}{|f|^p}\Bigg)^{\frac{q}{p}}$$
holds for some constant $C >0.$
\end{itemize}
\end{theorem}
Next we define an weight class similar to that of Sawyer's testing type condition. We use this testing condition to obtain the weak $(p, q)$ estimate for the minimal function.
\begin{Definition}\label{fmd4}
Given $0 < p \leq q < \infty,$ we say that the pair $(U, V) \in (W_{p, q}^+)^*$ if for each interval $I$ in $\mathbb{R}$ the following inequality
$$\int_{I} \frac{U}{m_{\mu}^+(\omega/\chi_I)^q} \leq C \Bigg( \int_{I}\omega \Bigg)^{\frac{q}{p}}$$
holds for some constant $C >0.$
\end{Definition}
The two weight strong $(p, q)$ estimate for the fractional one-sided minimal function is obtained using the $(W_{p, q}^+)^*$ condition and the main result is the following.

\begin{theorem}\label{fmt3}
Given $0 < p \leq q < \infty$ and for a pair of weight $(U, V)$ on $\mathbb{R}$ the following statements are equivalent.
\begin{itemize}
\item[(i)] The pair $(U,V) \in (W_{p, q}^+)^*.$
\item[(ii)] The strong $(p, q)$ holds for the fractional one-sided minimal function. \textit{i.e.}
$$ \int_{\mathbb{R}} \frac{U}{(m_{\mu}^+(f))^q} \leq C \Bigg( \int_{\mathbb{R}} \frac{V}{|f|^p} \Bigg)^{\frac{q}{p}}$$
holds for some constant $ C > 0.$
\end{itemize}
\end{theorem}
The next result is about the equivalence relationship between the weight class $W_{p, q}^+$ and $(W_{p, q}^+)^*.$ This result is obtained  as a corollary from the next two theorems. Next we provide these two theorems and finally we present the result regarding the equivalence of the weight class. 
\begin{theorem}\label{fmt4}
We assume that the pair $(U, V) \in W_{p,q}^+, 0 < p\leq q < \infty.$ Let $I=[\alpha, \beta]$ be any interval in $\mathbb{R}.$ Then
$$\int_{I^-} \frac{U}{m_{\mu}^+(\omega/ \chi_{I})^q} \leq C \Bigg( \int_{I^- \cup I^+} \omega\Bigg)^{\frac{q}{p}}$$
holds for $I^-$ and $I^+,$ where $I^-=[\alpha, \gamma]$ and $I^+=[\gamma, \delta]$ with $2|I^-|=|I|= 4|I^+|.$ 
\end{theorem}
\begin{theorem}\label{fmt5}
Let $\{K_l\}_{l \geq 0}$ be a decreasing sequence of nested intervals with the property that $|K_l| \to 0$ as $l \to \infty.$ Then for the pair $(U, V) \in W_{p, q}^+, 0< p \leq q < \infty,$
$$\lim_{l\to \infty} \int_{K_l} \frac{U}{m_{\mu}^+(\omega/ \chi_{K_l})^q}=0$$
holds.
\end{theorem}
\begin{corollary}\label{fmc1}
Let $0 < p \leq q < \infty.$ Then the weight classes $W_{p,q}^+$ and $(W_{p,q}^+)^*$ are equivalent.
 \textit{i.e.} a pair $(U, V) \in W_{p,q}^+$ if and only if $(U, V) \in (W_{p,q}^+)^*.$ 
\end{corollary}
\section{Proof of the Results}
\textbf{Theorem \ref{fmt1}.}
\begin{proof}
Let $(U, V) \in W_{p, q, \eta}^+, 0 < \eta <1.$ We choose $\eta= \frac{1}{2},$ then by the $W_{p, q, \frac{1}{2}}^{+}$ condition,
$$\frac{1}{|I^-|}\int_{I^-} U \leq C 2^{1 +(1 + \mu)q}\frac{1}{ |I|^{1 + (\mu - \frac{1}{p})q}} \Bigg(\frac{1}{|I|} \int_{I} V^{\frac{1}{p+1}}  \Bigg)^{\frac{(p+1)q}{p}},$$
where $I=[\alpha, \beta], I^-=[\alpha, \gamma]$ and $|I|=2|I^-|.$ Hence the pair $(U, V) \in W_{p, q}^+.$

Conversely, we assume that $(U, V) \in W_{p,q}^{+}.$ We need to show that $(U, V) \in W_{p, q, \eta}^+, 0< \eta < 1.$ We break the proof into two parts.\\
\textbf{Case I.} When $0 < \eta= \frac{|I^-|}{|I|} < \frac{1}{2},$ where $I=[\alpha, \beta]$ and $I^-=[\alpha, \gamma].$ 
Suppose $J=[\alpha, \delta]$ be a subinterval of $I$ with $2|J|=|I|.$ By the $W_{p, q}^+$ condition,
$$\int_{J}U\leq \frac{C|J|}{|I|^{1 +(\mu - \frac{1}{p})q}} \Bigg(\frac{1}{|I|}\int_{I} V^{\frac{1}{p+1}} \Bigg)^{\frac{(p+1)q}{p}}.$$
Since $2|J|=|I| > 2|I^-|.$ Thus $I^- \subset J.$ From the above inequality, we have
\begin{align*}
\frac{1}{|I^-|}\int_{I^-}U \leq \frac{1}{|I^-|}\int_{J}U &\leq \frac{C|J|}{|I^-||I|^{1 +(\mu - \frac{1}{p})q}} \Bigg(\frac{1}{|I|}\int_{I} V^{\frac{1}{p+1}} \Bigg)^{\frac{(p+1)q}{p}}\\
& \leq \frac{C}{\eta(1 - \eta)^{(1 + \mu)q}}\frac{1}{|I|^{1 +(\mu - \frac{1}{p})q}}\Bigg( \frac{1}{|I|}\int_{I}V^{\frac{1}{p+1}}\Bigg)^{\frac{(p+1)q}{p}}.
\end{align*}
Thus we get $(U, V)\in W_{p, q, \eta}^+, 0 < \eta \leq \frac{1}{2}.$\\
\textbf{Case II.} When $\frac{1}{2} < \eta =\frac{|I^-|}{|I|} < 1.$ We choose the smallest 
$N_0 \in \mathbb{N}$ such that $$\eta \leq 1 - \frac{1}{2^{N_0}}.$$ As $N_0$ is the smallest, so $$1 - \eta < \frac{1}{2^{N_0-1}}.$$ Let $\{J_k^-\}_k$ be a collection of subinterval of $I$ formed from the plus-minus decomposition of the interval $I.$ We get 
$$I^- \subset \cup_{k=1}^{N_0} J_k^-.$$
From the $W_{p, q}^+$ condition,
\begin{align*}
\int_{J_k^-}U &\leq \frac{C|J_k^-|}{(2|J_k^-|)^{1 + (\mu - \frac{1}{p})q}}\Bigg( \frac{1}{2|J_k^-|}\int_{I}V^{\frac{1}{p+1}} \Bigg)^{\frac{(p+1)q}{p}}\\
&= \frac{C|I||J_k^-||I|^{(1+\mu)q}}{(2|J_k^-|)^{1 + (\mu - \frac{1}{p})q}} \Bigg[ \frac{1}{|I|^{1+(\mu - \frac{1}{p})q}} \bigg( \frac{1}{|I|}\int_{I}V^{\frac{1}{p+1}}\bigg)^{\frac{(p+1)q}{p}}\Bigg]\\
&\leq C|I^-| \Bigg(\frac{I}{2|J_k^-|}\Bigg)^{(1 + \mu)q}\Bigg[ \frac{1}{|I|^{1+(\mu - \frac{1}{p})q}} \bigg( \frac{1}{|I|}\int_{I}V^{\frac{1}{p+1}}\bigg)^{\frac{(p+1)q}{p}}\Bigg].
\end{align*}
For $k \geq 1,$ we have
\begin{align*}
|J_k^-|&=|I|\bigg(1 - \frac{1}{2^k}\bigg)\\
\implies \frac{|I|}{2|J_k^-|}&=\frac{2^{k-1}}{2^k-1} \leq 2^{k-1}.
\end{align*}
Now
\begin{align*}
\frac{1}{|I^-|}\int_{I^-} U &\leq \frac{1}{|I^-|}\sum_{k=1}^{N_0} \int_{J_k^-}U\\
& \leq C \sum_{k=1}^{N_0}2^{(k-1)(1+\mu)q}\Bigg[ \frac{1}{|I|^{1+(\mu - \frac{1}{p})q}} \bigg( \frac{1}{|I|}\int_{I}V^{\frac{1}{p+1}}\bigg)^{\frac{(p+1)q}{p}}\Bigg]\\
& \leq \frac{CN_0}{2}\frac{1}{\eta(1 - \eta)^{(1 + \mu)q}}\Bigg[ \frac{1}{|I|^{1+(\mu - \frac{1}{p})q}} \bigg( \frac{1}{|I|}\int_{I}V^{\frac{1}{p+1}}\bigg)^{\frac{(p+1)q}{p}}\Bigg].
\end{align*}
Thus $(U, V) \in W_{p, q, \eta}^+, \frac{1}{2} < \eta < 1.$\\
This concludes the proof.
\end{proof}
\textbf{Theorem \ref{fmt2}}
\begin{proof}
$(i) \implies (ii).$\\
It is sufficient to prove the result for the function $f$ with $\frac{1}{f}$ having compact support\cite{CRUZ2}. We assume that $(U, V) \in W_{p, q}^+, 0 < p \leq q < \infty.$ For each $\lambda > 0,$ we consider the set
$$O_{\lambda}= \{x \in \mathbb{R}: m_{\mu}^+(f)(x) < \frac{1}{\lambda}\}.$$
Thus there exists a disjoint sequence of bounded and open interval $\{ I_k\}_{k \geq 1}$ such that
$$O_{\lambda} = \cup_{k \geq 1}I_k$$
and for the plus-minus decomposition of $I_k$
$$\frac{1}{|J_l^+|^{1 +  \mu}} \int_{J_l}|f| \leq \frac{8^{1+\mu}}{\lambda}.$$
Now, for a fixed $k,$
\begin{align*}
\int_{I_k}U &= \sum_{l}\int_{J_l^-}U\\
& \leq \frac{3^{(1+\mu)q}C}{{\lambda}^q} \sum_{l}\frac{1}{|J_l|^{(1+ \mu)q}} \Bigg(\int_{J_l} V^{\frac{1}{p+1}}\Bigg)^{\frac{(p+1)q}{p}}\Bigg(\frac{(8|J_l^+|)^{(1+\mu)q}}{\int_{J_l}|f|}\Bigg)^{q}\\
& \leq \frac{C}{{\lambda}^q} \sum_{l}\Bigg(\frac{3|J_l^+|}{|J_l|}\Bigg)^{(1+\mu)q} \Bigg(\int_{J_l} \frac{V^{\frac{1}{p+1}}}{|f|^{\frac{p}{p+1}}} |f|^{\frac{p}{p+1}}\Bigg)^{\frac{(p+1)q}{p}}\Bigg(\int_{J_l}|f|\Bigg)^{-q}\\
& \leq \frac{C}{{\lambda}^q} \sum_{l}\Bigg( \int_{J_l} \frac{V}{|f|^p}\Bigg)^{\frac{q}{p}}\Bigg(\int_{J_l}|f|\Bigg)^{q}\Bigg(\int_{J_l}|f|\Bigg)^{-q}\\
&=  \frac{C}{{\lambda}^q} \sum_{l} \Bigg( \int_{J_l} \frac{V}{|f|^p}\Bigg)^{\frac{q}{p}}\\
& \leq \frac{C}{{\lambda}^q}  \Bigg( \sum_{l} \int_{J_l} \frac{V}{|f|^p}\Bigg)^{\frac{q}{p}}\\
& \leq \frac{C}{{\lambda}^q}  \Bigg( \int_{I_k} \frac{V}{|f|^p}\Bigg)^{\frac{q}{p}}.
\end{align*}
Therefore
\begin{align*}
U(O_{\lambda})&=U(\cup_{k}I_k)\\
&= \sum_{k}U(I_k)\\
&\leq \sum_{k} \frac{C}{{\lambda}^q}  \Bigg( \int_{I_k} \frac{V}{|f|^p}\Bigg)^{\frac{q}{p}}\\
& \leq \frac{C}{{\lambda}^q} \Bigg( \int_{\mathbb{R}} \frac{V}{|f|^p}\Bigg)^{\frac{q}{p}}.
\end{align*}
$(ii) \implies (i).$\\
Let $I=[\alpha, \beta]$ be a finite interval in $\mathbb{R}.$ Let $I^-=[\alpha, \gamma]$ and $I^+=[\gamma, \beta]$ be two subintervals of $I$ with $|I^-|=|I^+|.$ If $x \in I^-,$ then
$$m_{\mu}^+(f)(x) \leq \frac{1}{|I^+|^{1+\mu}}\int_{I}V^{\frac{1}{p+1}},$$
where we assume $f = V^{\frac{1}{p+1}}/\chi_{I}.$\\
We choose $\lambda> 0$ such that
$$ \frac{1}{\lambda}=\frac{1}{|I^+|^{1 +\mu}}\int_{I}V^{\frac{1}{p+1}}.$$
From the assumption, we have
\begin{align*}
U(I^-) &\leq U\Bigg(\bigg\{ x \in \mathbb{R}: m_{\mu}^+(f)(x) < \frac{1}{\lambda} \bigg\}\Bigg)\\
& \leq \frac{C}{{\lambda}^q}\Bigg(\int_{I} \frac{V}{V^{\frac{p}{p+1}}} \Bigg)^{\frac{q}{p}}\\
&= C \Bigg(\frac{1}{|I^+|^{1+\mu}}\int_{I} V^{\frac{1}{p+1}}\Bigg)^q \Bigg( \int_{I} V^{\frac{1}{p+1}}\Bigg)^{\frac{q}{p}}\\
&= C \frac{|I^-|}{|I|^{1+(\mu - \frac{1}{p})q}} \Bigg( \int_{I} V^{\frac{1}{p+1}}\Bigg)^{\frac{(p+1)q}{p}}.
\end{align*}
Thus the pair $(U, V) \in W_{p, q}^+.$
\end{proof}

\textbf{Theorem \ref{fmt3}}
\begin{proof}
$(i) \implies (ii).$\\
As similar to the previous theorem, we only prove the statement for the function $f$ with $\frac{1}{f}$ having compact support.

Let $(U,V) \in (W_{p,q}^+)^*, 0 < p \leq q < \infty.$ For each $k \in \mathbb{Z},$ we define
$$\Omega_k=\{x \in \mathbb{R}: m_{\mu}^+(f)(x) < \frac{1}{2^k}\}.$$
Then by the definition of $m_{\mu}^+,$ there exists a disjoint sequence of bounded open intervals $\{I_{j,k}\}_j$ such that
$$\Omega_k= \cup_{j}I_{j, k}$$
and for each $x \in I_{j, k}=(\alpha_{j, k}, \beta_{j, k}),$
$$\int_{x}^{\beta_{j, k}} |f| \leq 2^{1 + \mu -k}|\beta_{j, k} - x|^{1+\mu}.$$
We construct an another pairwise disjoint set $\Omega_{j, k}$ for each integer $j$ and $k$ as 
$$\Omega_{j,k}= \Bigg\{ x \in I_{j, k}: m_{\mu}^+(f)(x) \geq \frac{1}{2^{k+1}} \Bigg\}.$$
Then
\begin{align*}
\int_{\mathbb{R}}\frac{U}{(m_{\mu}^+(f))^q} &= \sum_{j, k} \int_{\Omega_{j, k}} \frac{U}{m_{\mu}^+(f)^q}\\
&\leq\sum_{j, k} \int_{\Omega_{j, k}}2^{(k+1)q}U\\
& \leq 2^{q(2 + \mu)} \sum_{j, k}\int_{\Omega_{j, k}}\Bigg( \frac{1}{|\beta_{j, k} -x|^{1+\mu}}\int_{x}^{\beta_{j, k}}|f|\Bigg)^{-q}U.
\end{align*}
Using the technique from \cite{CRUZ2}, we assume that $\tau$ be a measure defined on $X = \mathbb{Z}\times\mathbb{Z}\times \mathbb{R}$ by $\nu \times \nu \times m,$ where $\nu$ stands for the counting measure and $m$ is the Lebesgue measure. We define $\phi$ on $X$ as
$$\phi(j, k, x) = \Bigg( \frac{1}{|\beta_{j, k} - x|^{1+\mu}} \int_{x}^{\beta_{j, k}} \omega dy \Bigg)^{-q}\chi_{\Omega_{j, k}}(x)U(x).$$
We define the operator $G$ as
$$G(g)(j, k, x)= \chi_{\Omega_{j, k}}(x)\frac{1}{\int_{x}^{\beta_{j, k}}g \omega dy} \int_{x}^{\beta_{j, k}} \omega dy$$
and the operator $H$ as
$$H(g)(j, k, x)= \chi_{\Omega_{j, k}}(x)\frac{1}{\int_{x}^{\beta_{j, k}} \omega dy} \int_{x}^{\beta_{j, k}} g \omega dy.$$
For $s >1,$ we have
\begin{align*}
\frac{1}{\int_{x}^{\beta_{j, k}}g \omega dy} \int_{x}^{\beta_{j, k}} \omega dy &= \frac{1}{\int_{x}^{\beta_{j, k}}g \omega dy} \bigg(\int_{x}^{\beta_{j, k}} \omega dy \bigg)^{1 - s +s}\\
&= \frac{1}{\int_{x}^{\beta_{j, k}}g \omega dy} \bigg(\int_{x}^{\beta_{j, k}} \omega dy \bigg)^{1 - s} \bigg(\int_{x}^{\beta_{j, k}} g^{\frac{1}{s}} {\omega}^{\frac{1}{s}} g^{- \frac{1}{s}} {\omega}^{1- \frac{1}{s}}dy \bigg)^{s}\\
& \leq \frac{1}{\int_{x}^{\beta_{j, k}}g \omega dy} \bigg( \int_{x}^{\beta_{j, k}}g \omega dy \bigg)\Bigg( \frac{1}{\int_{x}^6{\beta_{j,k}}} \bigg( \int_{x}^{\beta_{j, k}}g^{1 - s^{\prime}}\omega dy\bigg)\Bigg)^{s -1}\\
&=\Bigg( \frac{1}{\int_{x}^6{\beta_{j,k}}} \bigg( \int_{x}^{\beta_{j, k}}g^{1 - s^{\prime}}\omega dy\bigg)\Bigg)^{s -1}.
\end{align*}
Thus
$$G(g)(j,k,x) \leq \big(H(g^{1 - s^{\prime}})(j, k, x)\big)^{s-1}.$$
We set $s=1 + \frac{q}{p^2}$ and assume that $H$ maps from $L^{\frac{q}{p}}(\omega)$ to $L^{\frac{q^2}{p^2}}(X, \phi d\tau).$ Then from the above inequality 
\begin{align*}
2^{q(2 +\mu)}\int_{X} G\bigg(\frac{|f|}{\omega}\bigg)^{q} \phi d\tau & \leq 2^{q(2 +\mu)}\int_{X}H\bigg( \Big(\frac{\omega}{|f|}\Big)^{s^{\prime} -1} \bigg)^{(s-1)q} \phi d\tau\\
&\leq 2^{q(2 +\mu)}\Bigg(\int_{\mathbb{R}}\bigg( \Big(\frac{\omega}{|f|}\Big)^{s^{\prime} -1} \bigg)^{\frac{q}{p}}\omega dx\Bigg)^{\frac{p}{q}\times \frac{q^2}{p^2}}\\
&\leq C \bigg( \int_{\mathbb{R}} \frac{V}{|f|^p}dx \bigg)^{\frac{q}{p}}.
\end{align*}
Hence the result follows.
Therefore it remains to show the boundedness of $H.$ Following the argument stated in \cite{SAW}, it is sufficient to show that $H$ is weak $(1, \frac{q}{p}).$ That is for each $\lambda > 0,$ we need to show that
$$\int_{\{|H(g)| > \lambda\}}\phi d\tau \leq C \bigg( \frac{1}{\lambda} \int_{}g\omega dx\bigg)^{\frac{q}{p}}.$$
We define $t_{jk}(\lambda) = \inf \Gamma_{j, k}(\lambda),$ where $\Gamma_{j, k}(\lambda)=\{ x \in \Omega_{j, k}: H(g)(j, k, x) > \lambda\}$ and $J_{jk}=J_{jk}(\lambda)=[t_{jk}(\lambda), \beta_{jk}).$ Then
$$\frac{1}{\int_{J_{jk}}\omega}\int_{J_{jk}}g \omega \geq \lambda$$
and the pair of set $J_{jk}$ and $J_{lm}$ are either disjoint or one contains the other.\\
Let $\{J_i\}$ be the maximal elements of the family $\{J_{jk}\}$ so that $J_i$'s are disjoint. Thus
\begin{align*}
\int_{\{|H(g)| > \lambda\}}\phi d\tau &= \sum_{k, j}\int_{\Gamma_{jk}(\lambda)}\Bigg( \frac{1}{|\beta_{j,k} -x|^{1+\mu}}\Bigg)^{-q}Udx\\
&= \sum_{i}\sum_{\{(k, j): J_{jk} \subset {J_i}\}}\int_{\Gamma_{jk}(\lambda)}\Bigg( \frac{1}{|\beta_{j,k} -x|^{1+\mu}}\Bigg)^{-q}Udx\\
& \leq \sum_{i}\int_{J_i}\frac{U}{m_{\mu}^+(\omega/\chi_{J_i})^q}dx\\
& \leq C \sum_{i}\bigg(\int_{J_i} \omega dx\bigg)^{\frac{q}{p}}\\
& \leq C \sum_{i}\bigg( \frac{1}{\lambda} \int_{J_i}g\omega dx\bigg)^{\frac{q}{p}}\\
&\leq C\bigg(\frac{1}{\lambda} \int_{\mathbb{R}} g \omega dx \bigg)^{\frac{q}{p}}.
\end{align*}
This proves $H$ is weak $(1, \frac{q}{p})$ and hence concludes the proof.\\
$(ii) \implies (i).$\\
Let $f = \omega/\chi_{I},$ where $I$ is a fix interval and $\omega = V^{\frac{1}{p+1}}.$ And hence the result follows immediately.
\end{proof}

\textbf{Theorem \ref{fmt4}}
\begin{proof}
For each $\lambda > 0,$ we define

$$\Omega_{\lambda}= \Big\{ x \in I^-: m_{\mu}^+(\omega/ \chi_I)(x) < \frac{1}{\lambda}\Big\}.$$
There exists a sequence of disjoint intervals $\{I_j\}$ such that $\Omega_{\lambda}= \cup_{j}I_j.$\\
Also we have,
\begin{align*}
\int_{I^-} \frac{U}{m_{\mu}^+(\omega/ \chi_I)^q}dx &= q \int_0^{\infty}{\lambda}^{q-1}U(\Omega_{\lambda})d{\lambda}\\
&= I_1 + I_2,
\end{align*}
where, $$I_1= q \int_{0}^{\epsilon}{\lambda}^{q-1}U(\Omega_{\lambda})d\lambda$$
and
$$I_2= q \int_{\epsilon}^{\infty}{\lambda}^{q-1}U(\Omega_{\lambda})d\lambda,$$
for some $\epsilon > 0$ and the value of $\epsilon$ to be chosen later.
Now, 
\begin{align}
I_1 &= q \int_{0}^{\epsilon}{\lambda}^{q-1}U(\Omega_{\lambda})d\lambda \nonumber\\
& \leq U(I^-)q \int_{0}^{\epsilon}{\lambda}^{q-1} \nonumber\\
& = U(I^-){\epsilon}^q. \label{fme1}
\end{align}
For a fixed interval $\{I_j\},$ we construct two sequences of intervals $\{J_k^-\}$ and $\{J_k^+\}$ as the plus-minus decomposition of the interval $I_j.$ Then
\begin{align}
U(I_j) &= \sum_{k=1}^{\infty}U(J_k^-)\nonumber\\
& \leq 3^{(1+\mu)q}C \sum_{k=1}^{\infty}\frac{3|J_k^-|}{2|J_k|^{1+(1+\mu)q}}\bigg(\int_{J_k}\omega \bigg)^{\frac{(p+1)q}{p}}\nonumber\\
& = C \sum_{k=1}^{\infty}\bigg(\frac{1}{|J_k|^{(1+\mu)}} \int_{J_k}\omega \bigg)^q \bigg(\int_{J_k}\omega \bigg)^{\frac{q}{p}}\nonumber\\
& \leq C \sum_{k=1}^{\infty} \bigg(\frac{1}{{\lambda}}\bigg)^{\frac{(p+1)q}{p}}|J_k|^{(1 + \mu)\frac{q}{p}}\nonumber\\
& \leq C  \bigg(\frac{1}{{\lambda}}\bigg)^{\frac{(p+1)q}{p}}|I_j|^{(1 + \mu)\frac{q}{p}}. \label{fme2}  
\end{align}
Thus we obtain the estimate for the integral $I_2$ as
\begin{align}
 I_2 &= q \int_{\epsilon}^{\infty}\lambda^{q-1}U(\Omega_{\lambda})d\lambda \nonumber\\
 &= q \int_{\epsilon}^{\infty}\lambda^{q-1}\sum_{j}U(I_j)d\lambda \nonumber\\
 & \leq  q \int_{\epsilon}^{\infty}\lambda^{q-1}\bigg( \frac{1}{\lambda}\bigg)^{\frac{(p+1)q}{p}}|I^-|^{(1+\mu)\frac{q}{p}} \nonumber\\
 & \leq C \bigg(\frac{1}{\epsilon}\bigg)^{\frac{q}{p}}|I^-|^{(1+\mu)\frac{q}{p}}.\label{fme3}
\end{align}
We choose $\epsilon$ such that
\begin{equation}
\epsilon^q = \frac{(\omega(I^- \cup I^+))^{\frac{q}{p}}}{U(I^-)}. \label{fme4}
\end{equation}
As the pair $(U, V) \in W_{p, q}^+$ 
\begin{align}
U(I^-) &\leq C \frac{3^{1+(1 +\mu)q}|I^-|}{2|I^- \cup I^+|^{1 + (1 +\mu)q}}\bigg( \omega(I^- \cup I^+)\bigg)^{\frac{(p+1)q}{p}} \nonumber\\
&\leq \frac{C}{|I^-|^{(1+\mu)q}}\bigg( \omega(I^- \cup I^+)\bigg)^{\frac{(p+1)q}{p}}. \label{fme5}
\end{align}
From the inequality (\ref{fme4}) and (\ref{fme5}), we obtain the following.
\begin{equation}
\bigg(\frac{1}{\epsilon}\bigg)^{\frac{q}{p}} \leq \frac{C}{|I^-|^{( 1 +\mu)\frac{q}{p}}}\bigg( \omega(I^- \cup I^+)\bigg)^{\frac{q}{p}}. \label{fme6}
\end{equation}
Using the inequality (\ref{fme1}), (\ref{fme3}) and (\ref{fme6}), we conclude that 
$$I_1, I_2 \leq (\omega(I^- \cup I^+))^{\frac{q}{p}}.$$
It completes the proof.
\end{proof}

\textbf{Theorem \ref{fmt5}}
\begin{proof}
As the sequence $\{K_l\}$ is decreasing, so, $K_{l+1} \subset K_l, l \geq 0.$  Then for $0 \leq \mu < \infty,$ the minimal functions satisfies 
$$m_{\mu}^+(\omega/\chi_{K_l}) \leq m_{\mu}^+(\omega/ \chi_{K_{l+1}}), ~ k\geq 0.$$
We have$$1/ m_{\mu}^+(\omega/\chi_{K_l}) \to 0$$
\textit{a.e.} on $K_0$ as $l\to \infty.$\\
We set $I^-=K_0$ and by the the Theorem \ref{fmt4}, we obtain
$$\int_{K_0}\frac{U}{m_{\mu}^+(\omega/ \chi_{K_0})} < \infty.$$
Using the decreasing property of the sequence of interval $\{K_l\}$ and from the dominated convergence theorem,
$$\lim_{l\to \infty}\int_{K_l} \frac{U}{m_{\mu}^+(\omega/ \chi_{K_l})} \leq \lim_{l\to \infty}\int_{K_0} \frac{U}{m_{\mu}^+(\omega/ \chi_{K_l})}=0.$$ 
Hence the proof completes.
\end{proof}
\textbf{Corollary \ref{fmc1}}
\begin{proof}
Clearly $(W_{p, q}^+)^*$ implies $W_{p, q}^+.$ We only need to prove that $W_{p, q}^+$ implies $(W_{p, q}^+)^*.$ 

Let $(U, V)\in W_{p, q}^+.$ We show that for any interval $I=[\alpha, \beta]$ on $\mathbb{R},$ the inequality
$$\int_{I}\frac{U}{m_{\mu}^+(\omega/ \chi_{I})^q} \leq C  \bigg(\omega(I)\bigg)^{\frac{q}{p}}$$
holds for some constant $C >0.$

We construct a sequence of intervals $\{J_k^{\prime}\}$ as, for each $k \geq 1,$ we define $J_{k}^{\prime}=[x_k, \beta],$ where $\{x_k\}$ is the sequence constructed to form the plus-minus decomposition, $J_k^-$ and $J_k^+,$ of the interval $I$[see the Definition \ref{fmd3}]. Then $I =J_1^- \cup J_1^{\prime}$ and thus the following inequality follows.
$$\int_{I}\frac{U}{m_{\mu}^+(\omega/ \chi_{I})^q} \leq \int_{J_1^-} + \int_{J_1^{\prime}}.$$
From the Theorem \ref{fmt4},$$\int_{J_1^-} \leq C \bigg( \omega(J_1)\bigg)^{\frac{q}{p}}.$$
For the second integral, we can write $J_1^{\prime}= J_2^- \cup J_2^{\prime}$ and we repeat the process.
Continuing this process up to $l$ many times, we get,
$$\int_{I}\frac{U}{m_{\mu}^+(\omega/ \chi_{I})^q} \leq C \sum_{k=1}^{l}\bigg( \omega(J_k)\bigg)^{\frac{q}{p}} + \int_{J_l^{\prime}}\frac{U}{m_{\mu}^+(\omega/ \chi_{J_{l}^{\prime}})^q}$$
Using the Theorem \ref{fmt5}, we conclude that the last term tends to $0$ as $l \to \infty.$ Now, letting $l \to \infty,$
$$\int_{I}\frac{U}{m_{\mu}^+(\omega/ \chi_{I})^q} \leq C \bigg(\sum_{k=1}^{\infty} \omega(J_k)\bigg)^{\frac{q}{p}} \leq C \bigg( \omega(I)\bigg)^{\frac{q}{p}}.$$
This completes the proof.
\end{proof}

\subsection*{Acknowledgements}
D. Chutia was supported by the DST INSPIRE  
(Grant No. DST/INSPIRE Fellowship/2017/IF170509).
R. Haloi was supported by the DST MATRICS 
(Grant No. SERB/F/12082/2018-2019).

\end{document}